\documentclass[12]{article}
\usepackage{amsfonts}
\usepackage{amsmath}
\usepackage[unicode,colorlinks,plainpages=false]{hyperref}
\usepackage{theorem}
\newtheorem{example}{Example}[section]
\newtheorem{defn}{Definition}[section]
\newtheorem{thm}{Theorem}[section]

\newtheorem{lem}{Lemma}[section]

\newtheorem{proposition}{Proposition}[section]

\newcommand{\openbox}{$\begin{array}{c}
\hspace*{-0.55em}\sqcap \hspace*{-0.60em}\\[-0.4em] \hline
\multicolumn{1}{c}{\hspace*{-0.60em}}\\[-0.8em]
\end{array}
$}

\usepackage{enumerate}
\begin{document}

\centerline{\bf LEFT EQUALIZER SIMPLE SEMIGROUPS \footnote{ The paper will be published in Acta Mathematica Hungarica}}

\medskip

\centerline{ATTILA NAGY}
\centerline{Department of Algebra, Mathematical Institute}
\centerline{Budapest University of Technology and Economics}
\centerline{1521 Budapest, PO Box 91}
\centerline{e-mail: nagyat@math.bme.hu}

\begin{abstract}
In this paper we characterize and construct semigroups whose right regular representation is a left cancellative semigroup. These semigroups will be called
left equalizer simple semigroups. For a congruence $\varrho$ on a semigroup $S$,
let ${\mathbb F}[\varrho]$ denote the ideal of
the semigroup algebra ${\mathbb F}[S]$ which determines the kernel of the extended homomorphism of ${\mathbb F}[S]$ onto
${\mathbb F}[S/\varrho]$ induced by the canonical homomorphism of $S$ onto $S/\varrho$. We examine the right colons
$({\mathbb F}[\varrho]:_r{\mathbb F}[S])=\{ a\in {\mathbb F}[S]:\ {\mathbb F}[S]a\subseteq {\mathbb F}[\varrho ]\}$ in general, and in that special case
when $\varrho$ has the property that the factor semigroup $S/\varrho$ is left equalizer simple.
\end{abstract}

\section{Introduction and motivation}

Let $S$ be a semigroup. For an arbitrary element $a$ of $S$, let $\varrho _a$ denote the transformation of $S$ defined by $\varrho _a: s \mapsto sa$ ($s\in S$).
It is well known that $\theta _S=\{ (a,b)\in S\times S:\ (\forall x\in S)\ xa=xb\}$ is a congruence on $S$; this congruence is the kernel of the homomorphism
$\varphi : a \mapsto \varrho _a$ of $S$ into the semigroup of all right translations of $S$. The homomorphism $\varphi$ is called the right regular
representation of $S$; this is faithful if and only if $S$ is a left reductive semigroup (that is, whenever $xa=xb$ for some $a, b\in S$ and for all $x\in S$
then $a=b$). For convenience (as in \cite{Cr1:sg-3} or \cite{Cr2:sg-4}), the semigroup $\varphi (S)$ is also called the right regular representation of $S$.

It is an interesting problem to find couples $({\it C}_1, {\it C}_2)$ of classes ${\it C}_1$ and ${\it C}_2$ of semigroups for which the next assertion is true.
A semigroup $S$ belongs to ${\it C}_1$ if and only if the right regular representation of $S$ belongs to ${\it C}_2$.

It is easy to see that a semigroup $S$ is right commutative (that is, it satisfies the identity $xab=xba$ (\cite{Nagybook:sg-5})) if and only if the right
regular representation of $S$ is a commutative semigroup. A semigroup $S$ is medial (that is, it satisfies the identity $xaby=xbay$ (\cite{Nagybook:sg-5})) if
 and only if the right regular representation of $S$ is a left commutative semigroup (that is, it satisfies the identity $abx=bax$ (\cite{Nagybook:sg-5})).

In \cite{Y:sg-10}, a semigroup $S$ is called an $M$-inversive semigroup if, for each $a\in S$, there are elements $x, y\in S$ such that $ax$ and $ya$ are middle
units of $S$, that is, $caxd=cd$ and $cyad=cd$ is satisfied for all $c, d\in S$. In \cite{Cr2:sg-4}, it is proved that a semigroup $S$ is $M$-inversive if and
only if the right regular representation of $S$ is a right group (that is, a direct product of a right zero semigroup (satisfying the identity $ab=b$) and a group).

In Section 2, we define a special class of semigroups which contains all of the $M$-inversive ones. The semigroups belonging to this class will be called left
equalizer simple semigroups. Extending the above mentioned result of \cite{Cr2:sg-4}, we show that a semigroup $S$ is left equalizer simple if and only if the
right regular representation of $S$ is a left cancellative semigroup (that is, $xa=xb$ for some $x, a, b\in S$ implies $a=b$). We give a construction, and show
that a semigroup is left equalizer simple if and only if it is isomorphic to a semigroup defined by this construction.

Let $\it P$ be a property of semigroups. A congruence $\varrho$ on a semigroup $S$ will be called a $\it P$ congruence if the factor semigroup $S/\varrho$ has
the property $\it P$.

In Section 3, we examine left equalizer simple congruences on semigroups.
Let $S$ be a semigroup and $\varrho$ a congruence on $S$. As in \cite{Nagyleft:sg-6}, let $\varrho ^*$ denote the following congruence on $S$:
$(a, b)\in \varrho ^*$ for $a, b\in S$ if and only if $(xa, xb)\in \varrho$ for all $x\in S$.
By Corollary 1 of \cite{Nagyleft:sg-6}, \[(\cdot )\eta : \varrho \to \varrho ^*\] is a $\wedge$-homomorphism of ${\cal L}(S)$ into itself.
In Theorem~\ref{equiffcanc}, we show that $\varrho$ is a left equalizer simple congruence on a semigroup $S$ if and only if the congruence $\varrho ^*$ is a left cancellative
congruence on $S$.

Let $S$ be a semigroup and $\mathbb F$ a field. For an arbitrary congruence $\varrho$ on $S$, let ${\mathbb F}[\varrho]$ denote the ideal of
the semigroup algebra ${\mathbb F}[S]$ which determines the kernel of the extended homomorphism of the semigroup algebra ${\mathbb F}[S]$ onto the
semigroup algebra ${\mathbb F}[S/\varrho ]$  defined by the canonical homomorphism of $S$ onto the factor semigroup $S/\varrho$ (see \cite{Clifford1:sg-1}).
For an ideal $J$ of ${\mathbb F}[S]$, let $J|_S$ denote the restriction of the congruence on ${\mathbb F}[S]$ defined by $J$ to $S$. By Lemma 5 of Chapter 4
of \cite{Okninski:sg-9}, for every semigroup $S$ and every field $\mathbb F$, the mapping $J \mapsto \ J|_S$ is a surjective homomorphism
of the $\wedge$-semilattice $Id({\mathbb F}[S])$ of all ideals of ${\mathbb F}[S]$ onto the $\wedge$-semilattice of all congruences on $S$
such that ${\mathbb F}[\varrho ]|_S=\varrho$ for every congruence $\varrho$ on $S$. Let \[(\cdot )\kappa : \varrho \to {\mathbb F}[\varrho ].\]
By the above, $\kappa$ is an injective mapping of ${\cal L}(S)$ into $Id({\mathbb F}[S])$.

Using the notation of Section 3.6 of \cite{JespersOkn:sg-12} for semigroup algebras, if $J$ is an arbitrary ideal of ${\mathbb F}[S]$ then let
$(J:_r{\mathbb F}[S])=\{ a\in {\mathbb F}[S]: {\mathbb F}[S]a\subseteq J\}$. It is easy to see that $(J:_r{\mathbb F}[S])$ is an ideal of the semigroup
algebra ${\mathbb F}[S]$ such that $J\subseteq (J:_r{\mathbb F}[S])$. An ideal $(J:_r{\mathbb F}[S])$ will be called a right colon (more precisely, the right
colon of $J$ with respect to ${\mathbb F}[S]$).
It is a matter of checking to see that $(J_1\cap J_2:_r{\mathbb F}[S])=(J_1:_r{\mathbb F}[S])\cap (J_2:_r{\mathbb F}[S])$ for arbitrary ideals $J_1$ and $J_2$
of ${\mathbb F}[S]$. Thus
\[(\cdot )\Phi : J\to (J:_r{\mathbb F}[S])\] is a $\wedge$-homomorphism of $Id({\mathbb F}[S])$ into itself.

By the above, we can consider the following diagram.

\begin{equation}\label{commutative}
\begin{matrix}{\cal L}(S)&\overrightarrow{\empty \eta \empty}&{\cal L}(S)\\
\kappa \downarrow &\empty&\downarrow \kappa&\\
Id({\mathbb F}[S])&\overrightarrow{\empty \Phi \empty} &Id({\mathbb F}(S])\end{matrix}
\end{equation}

\medskip

We shall say that the diagram (1) is commutative for some congruence $\varrho$ on a semigroup $S$ if
\[(\varrho )(\eta \circ \kappa )=(\varrho )(\kappa \circ \Phi),\] that is,
 \[{\mathbb F}[\varrho ^*]=({\mathbb F}[\varrho ]:_r{\mathbb F}[S]).\] In Theorem~\ref{th4}, we prove that, for arbitrary left equalizer simple congruence on a
 semigroup $S$, the diagram (1) is commutative.

In Section 4, our results will be applied for the matrix representation of finite left equalizer simple semigroups. We also prove a theorem for arbitrary
(not necessarily finite) semigroups. We show that if $\varrho$ is a left equalizer simple congruence on a semigroup $S$ then the right colon
$({\mathbb F}[\varrho ]:_r{\mathbb F}[S])$ equals the augmentation ideal ${\mathbb F}[\omega _S]$ of the semigroup algebra ${\mathbb F}[S]$ if and only if the
factor semigroup $S/\varrho$ is an ideal extension of a left zero semigroup by a null semigroup.

\medskip

\section{Left equalizer simple semigroups}

Let $S$ be a semigroup and $H$ a non-empty subset of $S$. By the left equalizer of $H$ we mean the set of all elements $x$ of $S$ for which $|xH|=1$, that is,
$xa=xb$ is satisfied for all $a, b\in H$. It is clear that the left equalizer of $H$ is either empty or a left ideal of $S$.

\begin{lem}\label{lm2mind} On an arbitrary semigroup $S$, the following assertions are equivalent.
\begin{enumerate}
\item [(i)] The left equalizer of any two-element subset of $S$ is either empty or $S$.
\item [(ii)] The left equalizer of any subset of $S$ is either empty or $S$.
\end{enumerate}
\end{lem}

{\bf Proof}. Assume $(i)$. Let $H$ be a non-empty subset of $S$. If $|H|=1$ then the left equalizer of $H$ is $S$. Consider the case when $|H|\geq 2$.
If $x_0\in S$ is in the left equalizer of $H$ then, for every two elements $a, b\in H$ with $a\neq b$, we have $x_0a=x_0b$ and so $x_0$ belongs to the left
equalizer of the subset $\{ a, b\}$. Thus every $x\in S$ is in the left equalizer of $\{ a, b\}$, that is, $xa=xb$ for all $x\in S$.
As $a, b\in H$ are arbitrary, we have $|xH|=1$ for all $x\in S$ and so the left equalizer of $H$ is $S$.

It is obvious that condition $(ii)$ implies condition $(i)$.
\hfill\openbox

\begin{defn}\label{dfalmost} A semigroup $S$ will be called a left equalizer simple semigroup if, for arbitrary non-empty subset $H$ of $S$, the left
equalizer of $H$ is either empty or equal to $S$. Equivalently (see Lemma~\ref{lm2mind}), for arbitrary elements $a, b\in S$, the assumption $x_0a=x_0b$ for
some $x_0\in S$ implies $xa=xb$ for all $x\in S$.
\end{defn}

\medskip
A semigroup $S$ is said to be left simple if $S$ is the only left ideal of $S$. It is obvious that every left simple semigroup is left equalizer simple.

\begin{lem}\label{lmMinversive} Every $M$-inversive semigroup is left equalizer simple.
\end{lem}

{\bf Proof}. Let $x_0, a, b$ be arbitrary elements of an $M$-inversive semigroup $S$ with $x_0a=x_0b$. Then there is an element $y\in S$ such that $yx_0$ is a
middle unit of $S$ and so, for all $x\in S$, we have $xa=xyx_0a=xyx_0b=xb$. Hence $S$ is left equalizer simple.\hfill\openbox

\medskip

The next theorem is an extension of Theorem 1 of \cite{Cr2:sg-4}.

\begin{thm}\label{threpr} A semigroup $S$ is left equalizer simple if and only if the right regular representation of $S$ is a left cancellative semigroup.
\end{thm}

{\bf Proof}. Let $S$ be a left equalizer simple semigroup. Let $\varphi$ denote the canonical homomorphism of $S$ onto the right regular representation of $S$. If
\[\varphi (x_0)\varphi (a)=\varphi (x_0)\varphi (b)\] for some $x_0, a, b\in S$ then, for an arbitrary element $x_1\in S$, we have \[x_1x_0a=x_1x_0b.\]
As $S$ is left equalizer simple, it follows that $xa=xb$ for all $x\in S$. Hence $\varphi (a)=\varphi (b)$. Thus the right regular representation
of $S$ is a left cancellative semigroup.

Conversely, assume that the right regular representation $\varphi (S)$ of a semigroup $S$ is left cancellative. If $x_0a=x_0b$ for some $x_0, a, b\in S$
then \[\varphi (x_0)\varphi (a)=\varphi (x_0)\varphi (b)\] in $\varphi (S)$. As $\varphi (S)$ is left cancellative, we get
$\varphi (a)=\varphi (b)$ which means that $xa=xb$ for all $x\in S$. Thus $S$ is a left equalizer simple semigroup.
\hfill\openbox

\medskip

{\bf Construction 1}: Let $T$ be a left cancellative semigroup. For each $t\in T$, associate a nonempty set $S_t$ such that $S_t\cap S_r=\emptyset$ for every
$t\neq r$. As $T$ is left cancellative, $x \mapsto tx$ is an injective mapping of $T$ onto $tT$.

For arbitrary couple $(t, r)\in T\times T$ with $r\in tT$, let $f_{(t, r)}$ be a mapping of $S_t$ into $S_r$. For all $t\in T$, $r\in tT$, $q\in rT\subseteq tT$
and $a\in S_t$, assume \[(a)f_{(t,r)}\circ f_{(r,q)}=(a)f_{(t, q)}.\]

On the set $S=\cup _{t\in T}S_t$ define an operation $\star$ as follows: for arbitrary $a\in S_t$ and $b\in S_x$, let \[a\star b=(a)f_{(t, tx)}.\]

If $a\in S_t$, $b\in S_x$, $c\in S_y$ are arbitrary elements then
\[a\star (b\star c)=a\star (b)f_{(x, xy)}=(a)f_{(t, t(xy))}=\]
\[=(a)(f_{(t, tx)}\circ f_{tx, t(xy)}=(a)f_{(t, tx)}\star c=(a\star b)\star c.\] Thus
$(S; \star )$ is a semigroup.

Let $x\in T$ and $a\in S_x$ be arbitrary. Then, for arbitrary $b\in S_x$ and $c\in S_t$ ($t\in T$), we have \[c\star a=(c)f_{(t, tx)}=c\star b.\]
Thus $(a, b)\in \theta _S$ and so $S_x$ is contained by the $\theta _S$-class $[a]_{\theta _S}$ of $S$. If $d\in [a]_{\theta _S}\cap S_y$ ($y\in T$) then,
for an arbitrary $c\in S_t$ ($t\in T$), we have
\[(c)f_{(t, tx)}=c\star a=c\star d=(c)f_{(t, ty)}\] from which we get $tx=ty$. As $T$ is left cancellative, we get $x=y$ and so $d\in S_x$.
Consequently $S_x=[a]_{\theta _S}$. Hence each set $S_x$ ($x\in T$) is a $\theta _S$-class of $S$.

\medskip

\begin{example}\label{example1}\rm
Let $T=\{ 1, 2\}$ be a two-element group in which $1$ is the identity element. Let $S_1=\{ x_1, y_1\}$ and $S_2=\{ x_2, y_2\}$ be disjoint sets.
For arbitrary $t, s\in T$, let $f_{(t, ts)}$ defined by: $(x_t)f_{(t, ts)}=x_{ts}$ and $(y_t)f_{(t, ts)}=y_{ts}$. It is easy to see that the conditions of
Construction 1 are satisfied. The Cayley-table of the semigroup $(S; \star )$ is Table~\ref{Fig0}.
\begin{table}[htbp]
\begin{center}
\begin{tabular}{l|l l l l }
 &$x_1$&$y_1$&$x_2$&$y_2$\\ \hline
$x_1$&$x_1$&$x_1$&$x_2$&$x_2$\\
$y_1$&$y_1$&$y_1$&$y_2$&$y_2$\\
$x_2$&$x_2$&$x_2$&$x_1$&$x_1$\\
$y_2$&$y_2$&$y_2$&$y_1$&$y_1$\\
\end{tabular}
\caption{}\label{Fig0}
\end{center}
\end{table}
\end{example}

\begin{thm}\label{konstrukcio} A semigroup is a left equalizer simple semigroup if and only if it is isomorphic to a semigroup defined in Construction 1.
\end{thm}

\noindent
{\bf Proof}.  Let $(S; \star )$ be a semigroup defined in Construction 1. We show that $(S; \star )$ is left equalizer simple. Assume
\[x_0\star a=x_0\star b\] for some $x_0, a, b\in S$. Let $\xi, t, r\in T$ such element for which $x_0\in S_{\xi}$, $a\in S_t$ and $b\in S_r$.
Then \[(x_0)f_{ (\xi , \xi t)}=(x_0)f_{(\xi , \xi r)}\] and so \[\xi t=\xi r\] in $T$. As $T$ is left cancellative, $t=r$ and so $\eta t=\eta r$ for all
$\eta \in T$. Thus, for every $\eta \in T$ and $x\in S_{\eta }$,
\[x\star a=(x)f_{(\eta, \eta t)}=(x)f_{(\eta , \eta r)}=x\star b.\] Consequently $S$ is a left equalizer simple semigroup.

To prove the converse, let $S$ be a left equalizer simple semigroup. Then, by Theorem~\ref{threpr}, the right regular representation
$T=\varphi (S)=S/\theta _S$ is a left cancellative semigroup. Construct the semigroup as in Construction 1 defined by the semigroup $T=S/\theta _S$.
For arbitrary $t\in T$, let $S_t$ be the $\theta _S$-class of $S$ for which $\varphi (S_t)=t$. For arbitrary couple
$(t, r)\in T\times T$ with $r\in tT$, let $f_{(t, r)}$ be a mapping of $S_t$ into $S_r$ defined by the following way. As $T$ is left cancellative, there is a
unique element $x\in T$ such that $tx=r$. For arbitrary $a\in S_t$, let $(a)f_{(t, r)}=aw$, where $w$ is an arbitrary element of $S_x$.
For arbitrary $a\in S_t$ and $b\in S_x$, \[a\star b=(a)f_{(t, tx)}=ab.\] Hence, $S$ is isomorphic to the semigroup $(S; \star )$.
\hfill\openbox

\section{Left equalizer simple congruences}

According to our definition for the $\it P$ congruence from Section 1, a congruence $\varrho$ on a semigroup $S$ is a left equalizer simple
congruence if, for every $a, b\in S$, the assumption $(x_0a, x_0b)\in \varrho$ for some $x_0\in S$ implies $(xa, xb)\in \varrho$ for all $x\in S$.
A congruence $\varrho$ on a semigroup $S$ is a left cancellative congruence if, for every $a, b\in S$, the assumption
$(x_0a, x_0b)\in \varrho$ for some $x_0\in S$ implies $(a, b)\in \varrho$.

\begin{thm}\label{equiffcanc} A congruence $\varrho$ on a semigroup $S$ is a left equalizer simple congruence if and only if $\varrho ^*$ is a left
cancellative congruence on $S$.
\end{thm}

{\bf Proof}. Let $\varrho$ be an arbitrary congruence on a semigroup $S$. As $\varrho \subseteq \varrho ^*$, we can consider the congruence
$\varrho ^*/\varrho$ (see \cite{Howie:sg-11}). By Theorem 5.6 of \cite{Howie:sg-11}, $(S/\varrho )/(\varrho ^*/\varrho )\cong S/\varrho ^*$.
By Lemma 7 of \cite{Nagyleft:sg-6}, $\varrho ^*/\varrho = \iota ^*_{S/\varrho}$, where $\iota _{S/\varrho}$ denotes the identity relation on the factor
semigroup $S/\varrho$. As $\iota ^*_{S/\varrho}=\theta _{S/\varrho}$, we get that the semigroup $S/\varrho ^*$ is the right regular representation of
$S/\varrho$. Thus $\varrho$ is a left equalizer simple congruence on $S$ if and only if the factor semigroup $S/\varrho$ is a left equalizer simple semigroup.
By Theorem~\ref{threpr}, this last condition is equivalent to the condition that the right regular representation of $S/\varrho$, that is, the semigroup
$S/\varrho ^*$ is left cancellative. This last condition means that $\varrho ^*$ is a left cancellative congruence on $S$. \hfill\openbox

\bigskip

Let $S$ be a semigroup and $\mathbb F$ a field. Let ${\mathbb F}[S]$ denote the semigroup algebra of $S$ over $\mathbb F$. By page 159
of \cite{Clifford1:sg-1}, $S$ can be considered as the basis of ${\mathbb F}[S]$ and every element $a$ of ${\mathbb F}[S]$ can be written as a
sum $\sum _{s\in S}a(s)s$; this is a finite sum since only a finite number of coefficients $a(s)\in {\mathbb F}$ are $\neq 0$.
For arbitrary elements $a, b\in {\mathbb F}[S]$ and $\alpha \in {\mathbb F}$,
\[a+b=\sum _{s\in S}(a(s)+b(s))s,\]
\[ab=\sum _{s\in S}\left(\sum_{\{x,y\in S; xy=s\}}a(x)b(y)\right)s,\]
\[\alpha a=\sum _{s\in S}\left(\alpha (a(s))\right)s.\]

\medskip

It is easy to see that an element $a=\sum _{s\in S}a(s)s$ of ${\mathbb F}[S]$ is in  ${\mathbb F}[\varrho]$  if and only if
$\sum_{s\in A}a(s)=0$ for every $\varrho$-class $A$ of $S$. If $\omega _S$ denotes the universal relation on $S$ then the ideal ${\mathbb F}[\omega _S]$ is
the set of all elements $\sum_{s\in S}a(s)s$ of ${\mathbb F}[S]$ for which $\sum _{s\in S}a(s)=0$; this ideal is the (so called) augmentation ideal of
${\mathbb F}[S]$.

\begin{lem}\label{th2} For an arbitrary field $\mathbb F$ and an arbitrary congruence $\varrho$ on a semigroup $S$,
${\mathbb F}[\varrho ]\subseteq {\mathbb F}[\varrho ^*]\subseteq({\mathbb F}[\varrho ]:_r{\mathbb F}[S])$ (with the notations of the diagram (1),
$(\varrho )\kappa \subseteq (\varrho )(\eta \circ \kappa )\subseteq (\varrho )(\kappa \circ \Phi)$).
\end{lem}

{\bf Proof}. The algebra ideal ${\mathbb F}[\alpha ]$ of a semigroup congruence $\alpha$ is obviously generated (as an $\mathbb F$-space)
by differences $s-s'$, where $(s, s')\in \alpha$. Thus the inclusion ${\mathbb F}[\varrho ]\subseteq {\mathbb F}[\varrho ^*]$ is obvious.

If $(s, s')\in \varrho ^*$ then, for every $x\in S$, $(xs, xs')\in \varrho$, hence $xs-xs'\in {\mathbb F}[\varrho ]$. From this it follows that
${\mathbb F}[\varrho ^*]\subseteq({\mathbb F}[\varrho ]:_r{\mathbb F}[S])$.\hfill\openbox

\medskip

\begin{example}\label{example2}\rm
Let $S=\{a, b, c, d, e\}$ be a semigroup defined by Table~\ref{Fig1}
(see \cite{Bogdanovic:sg-2}; page 167, the last Cayley-table in row 7):

\begin{table}[htbp]
\begin{center}
\begin{tabular}{l|l l l l l }
 &$a$&$b$&$c$&$d$&$e$\\ \hline
$a$&$b$&$b$&$a$&$a$&$b$\\
$b$&$b$&$b$&$b$&$b$&$b$\\
$c$&$b$&$b$&$c$&$c$&$b$\\
$d$&$b$&$b$&$d$&$d$&$b$\\
$e$&$a$&$b$&$a$&$b$&$e$\\
\end{tabular}
\caption{}\label{Fig1}
\end{center}
\end{table}

As the columns of the Cayley-table are pairwise distinct, $S$ is left reductive and so the identity relation $\iota _S$ on $S$ is left reductive.
Thus, by Theorem 1 of \cite{Nagyleft:sg-6}, \[\iota _S=\iota ^*_S\] and so \[{\mathbb F}[\iota _S]=\{ 0\}={\mathbb F}[\iota ^*_S].\] It is a matter of
checking to see that
\[d+a-b-c\in ({\mathbb F}[\iota _S]:_r{\mathbb F}[S]).\] Thus
\[\{ 0\}\subset ({\mathbb F}[\iota _S]:_r{\mathbb F}[S]).\]
This example shows that ${\mathbb F}[\varrho]=({\mathbb F}[\varrho ]:_r {\mathbb F}[S])$ or
${\mathbb F}[\varrho ^* ]=({\mathbb F}[\varrho ]:_r {\mathbb F}[S])$ is not satisfied for an arbitrary congruence $\varrho$ on a semigroup $S$.
\end{example}
\bigskip

By Lemma~\ref{th2}, for an arbitrary field $\mathbb F$ and an arbitrary congruence $\varrho$ on a semigroup $S$:
\begin{itemize}
\item[(a)] ${\mathbb F}[\varrho ]\subseteq {\mathbb F}[\varrho ^*]$,
\item[(b)] ${\mathbb F}[\varrho ]\subseteq ({\mathbb F}[\varrho ]:_r {\mathbb F}[S])$,
\item[(c)] ${\mathbb F}[\varrho^* ]\subseteq ({\mathbb F}[\varrho ]:_r {\mathbb F}[S])$.
\end{itemize}

\medskip

In case $(a)$, the equation ${\mathbb F}[\varrho ]={\mathbb F}[\varrho ^*]$ holds if and only if $\varrho=\varrho ^*$, because the mapping
$\kappa : \varrho \to {\mathbb F}[\varrho ]$ is injective (by Lemma 5 of Chapter 4 of \cite{Okninski:sg-9}).
By Theorem 1 of \cite{Nagyleft:sg-6}, a congruence $\varrho$ on a semigroup $S$ is left reductive if and only if $\varrho =\varrho ^*$. Thus we have the
following proposition.

\begin{proposition} Let $\mathbb F$ be an arbitrary field. The equation ${\mathbb F}[\varrho ]={\mathbb F}[\varrho ^*]$ is satisfied for a
congruence $\varrho$ on a semigroup $S$ if and only if $\varrho$ is a left reductive congruence on $S$.\hfill\openbox
\end{proposition}

\medskip

To characterize congruences $\varrho$ on a semigroup $S$ for which the equation ${\mathbb F}[\varrho ]=({\mathbb F}[\varrho ]:_r {\mathbb F}[S])$ holds
(see case $(b)$), consider the notion of the right annihilator $Ann_r({\mathbb F}[S])$ of a semigroup algebra ${\mathbb F}[S]$.
Recall that $Ann_r({\mathbb F}[S])=\{a\in {\mathbb F}[S]: (\forall x\in {\mathbb F}[S])\ xa=0\}=(\{ 0\}:_r{\mathbb F}[S])$.
The right annihilator is said to be trivial if it contains only the zero $0$ of ${\mathbb F}[S]$. Here we refer to
\cite{Nagy0:sg-7} and \cite{NagyRonyai:sg-8}, where the finite semigroups $S$ with condition $Ann_r({\mathbb F}[S])=\{ 0\}$ were examined.

\begin{proposition}\label{th6} Let $\mathbb F$ be an arbitrary field. For a congruence $\varrho$ on a semigroup $S$, the equation
${\mathbb F}[\varrho ]=({\mathbb F}[\varrho ]:_r {\mathbb F}[S])$ is satisfied if and only if the right annihilator $Ann_r({\mathbb F}[S/\varrho])$ is
trivial.
\end{proposition}

{\bf Proof}.  Let $\tau$ denote the (extended) canonical homomorphism of ${\mathbb F}[S]$ onto ${\mathbb F}[S/\varrho]$. For an element
$a\in {\mathbb F}[S]$, $\tau (a)\in Ann_r({\mathbb F}[S/\varrho])$ if and only if, for all $x\in {\mathbb F}[S]$,
$\tau (xa)=\tau (x)\tau (a)=0$ which is equivalent to the condition that $xa\in {\mathbb F}[\varrho]$ for all $x\in {\mathbb F}[S]$.
This last condition means that $a\in ({\mathbb F}[\varrho ]:_r {\mathbb F}[S])$. From this it follows that
${\mathbb F}[\varrho ]=({\mathbb F}[\varrho ]:_r {\mathbb F}[S])$ if and only if $Ann_r({\mathbb F}[S/\varrho])$ is trivial.
\hfill\openbox

\medskip

The condition that equation ${\mathbb F}[\varrho ^*]=({\mathbb F}[\varrho ]:_r{\mathbb F}[S])$ holds in case (c) is equivalent to the condition that the
diagram (1) is commutative for $\varrho$. The next theorem shows that the left equalizer simplicity is a sufficient condition for a congruence $\varrho$ to
be the diagram (1) commutative for $\varrho$.

\begin{thm}\label{th4} Let $\mathbb F$ be an arbitrary field. If $\varrho$ is a left equalizer simple congruence on a semigroup $S$ then
${\mathbb F}[\varrho ^*]=({\mathbb F}[\varrho ]:_r{\mathbb F}[S])$, that is, $(\varrho )(\eta \circ \kappa )=(\varrho )(\kappa \circ \Phi)$ which means that
the diagram (1) is commutative for $\varrho$.
\end{thm}

{\bf Proof}.
Let $\varrho$ be a left equalizer simple congruence on a semigroup. We must show that
\[{\mathbb F}[\varrho ^*]=({\mathbb F}[\varrho ]:_r{\mathbb F}[S]).\]
By Lemma~\ref{th2}, it is sufficient to show that
\[({\mathbb F}[\varrho ]:_r{\mathbb F}[S])\subseteq {\mathbb F}[\varrho ^*].\] Let
\[\sum_{s\in S}a(s)s\in ({\mathbb F}[\varrho ]:_r{\mathbb F}[S])\] be arbitrary. Let $x\in S$ be an arbitrary element. Then
\[\sum_{s\in S}a(s)xs\in {\mathbb F}[\varrho ].\] Let $A$ be a $\varrho ^*$-class of $S$, and let $B$ denote the $\varrho$-class of $S$ containing $xA$.
We show that $xS\cap B=xA$. The inclusion $xA\subseteq xS\cap B$ is obvious. If $xs\in B$ for some $s\in S$ then $(xs, xa)\in \varrho$ for some $a\in A$.
As $\varrho$ is a left equalizer simple congruence on $S$, we get $(ts, ta)\in \varrho$ for every $t\in S$. Thus $(s, a)\in \varrho ^*$ and so $s\in A$.
This implies $xS\cap B\subseteq xA$. Consequently $xS\cap B=xA$. Thus the sum $\sum_{s\in A}a(s)$ is the same as the sum $\sum_{xs\in B}a(s)$.
The latter sum is $0$ by assumption, giving that the former sum is also $0$, which proves the theorem.\hfill\openbox

\medskip

\section{Matrix representations}

Let $S$ be a finite semigroup and $\mathbb F$ an arbitrary field.
By an {\it $S$-matrix} over $\mathbb F$ we mean a mapping of the
direct product $S\times S$ into $\mathbb F$.
The set ${\mathbb F}_{S\times S}$
of all $S$-matrices over $\mathbb F$
is an algebra over $\mathbb F$ under the usual addition and multiplication of matrices
and the product of matrices by scalars.

A homomorphism $\gamma$ of a $S$ into the multiplicative semigroup of the full matrix algebra ${\mathbb F}_{n\times n}$ of all $n\times n$ matrices over
$\mathbb F$ is called a matrix representation of $S$ of order $n$. We say that $\gamma$ is faithful if it is injective.

For an element $s\in S$,
let ${\bf R}^{(s)}$ denote the $S$-matrix defined by
\[ {\bf R}^{(s)}((x,y))=\begin{cases}
1, & \text{if $xs=y$}\\
0 & \text{otherwise,}
\end{cases}\]
where $1$ and $0$ denote the identity element and the zero element of $\mathbb F$, respectively.
This matrix is called the right matrix of the element $s$ of $S$.
The mapping
\[{\mathcal R}_{{\mathbb F}}:\ s\mapsto {\bf R}^{(s)}\]
is a matrix representation of $S$ over ${\mathbb F}$ (see, for example, Exercise 4 in
\S 3.5 of \cite{Clifford1:sg-1}). ${\mathcal R}_{{\mathbb F}}$ describes in terms of matrices ${\bf R}^{(s)}$ the maps
$x\mapsto xs$. Thus, it is essentially the right regular representation.
It is obvious that a semigroup containing an identity element is left reductive. Thus, for an arbitrary finite semigroup $S$,
the restriction ${\mathcal R'}_{\mathbb F}$ of the right regular matrix representation of $S^1$ to $S$ is a faithful matrix representation of $S$.
The matrices $\{ {\mathcal R'}_{\mathbb F}(s); s\in S\}$ are also linearly independent over $\mathbb F$.

If $S$ is an arbitrary finite $n$-element left reductive semigroup then the system $\{{\mathcal R}_{\mathbb F}(s); s\in S\}$
has $n$ pairwise different matrices, but these matrices are not linearly independent over $\mathbb F$, in general. By Definition 2.2 of \cite{NagyRonyai:sg-8},
a finite semigroup $S$ is called an ${\cal R}_{\mathbb F}$-independent semigroup if the system $\{{\mathcal R}_{\mathbb F}(s); s\in S\}$ of matrices is
linearly independent over $\mathbb F$.

For a congruence $\varrho$ on $S$, define the following sequence (see the diagram (\ref{commutative})):
\[\varrho ^{(n)}=(\varrho )\eta ^n,\ n=0, 1, \dots\]

If $\varrho$ is a left equalizer simple congruence on $S$ then $\varrho ^{(1)}$ is a left cancellative congruence on $S$ by Theorem~\ref{equiffcanc}.
As a left cancellative congruence is also  left reductive, we have
\[\varrho ^{(1)}=\varrho ^{(2)}=\cdots\] by Theorem 1 and Theorem 2 of \cite{Nagyleft:sg-6}.
As a left cancellative congruence is also left equalizer simple, we have
\[{\mathbb F}[\varrho ^{(2)}]=({\mathbb F}[\varrho ^{(1)}]:_r{\mathbb F}[S])\] by Theorem~\ref{th4}.
From this it follows that
\begin{equation}\label{2}
{\mathbb F}[\varrho ^{(1)}]=({\mathbb F}[\varrho ^{(1)}]:_r{\mathbb F}[S]).
\end{equation}
By Proposition~\ref{th6}, it means that $Ann_r({\mathbb F}[S/\varrho ^{(1)}])$ is trivial.
Then the factor semigroup $S/\varrho ^{(1)}$ is ${\cal R}_{\mathbb F}$-independent by Theorem 2.1 of \cite{NagyRonyai:sg-8}.

The above result is not too interesting if $S/\varrho ^{(1)}$ is a one-element semigroup, that is, if $\varrho ^{(1)}=\omega _S$ ($\omega _S$
denotes the universal relation on $S$).  This special case is equivalent to (see also (\ref{2}))
\[({\mathbb F}[\varrho ^{(1)}]:_r{\mathbb F}[S])={\mathbb F}[\varrho ^{(1)}]={\mathbb F}[\omega _S],\] where ${\mathbb F}[\omega _S]$ is the augmentation
ideal of ${\mathbb F}[S]$. In the last part of the paper we deal with the following question: What can we say about the factor semigroup $S/\varrho$
if $\varrho$ is a left equalizer simple congruence on a (not necessarily finite) semigroup $S$ such that the right colon
$({\mathbb F}[\varrho ]:_r{\mathbb F}[S])$ equals the augmentation ideal ${\mathbb F}[\omega _S]$ of ${\mathbb F}[S]$.

\begin{thm}\label{th5} Let $\varrho$ be a left equalizer simple congruence on a semigroup $S$ and $\mathbb F$ a field. The right colon
$({\mathbb F}[\varrho ]:_r{\mathbb F}[S])$ equals the augmentation ideal ${\mathbb F}[\omega _S]$ of ${\mathbb F}[S]$ if and only if the factor
semigroup $S/\varrho$ is an ideal extension of a left zero semigroup by a null semigroup (which means that $S$ has an ideal $L$ which is a left zero
semigroup and the Rees factor semigroup $S/L$ is a null semigroup (which means that $(S/L)^2=\{ 0\}$)).
\end{thm}

{\bf Proof}. As $\varrho$ is left equalizer simple,
${\mathbb F}[\varrho ^{(1)}]=({\mathbb F}[\varrho ]:_r{\mathbb F}[S])$ by Theorem~\ref{th4}. Thus
$({\mathbb F}[\varrho ]:_r{\mathbb F}[S])={\mathbb F}[\omega _S]$ is satisfied if and only if
${\mathbb F}[\varrho ^{(1)}]={\mathbb F}[\omega _S]$. By Lemma 5 of Chapter 4
of \cite{Okninski:sg-9}, this is equivalent to $\varrho^{(1)}=\omega _S$.
By the proofs of Theorem 5 and Theorem 6 of \cite{Nagyleft:sg-6}, this is equivalent to the condition that the factor semigroup $S/\varrho$ is an ideal
extension of a left zero semigroup by a null semigroup.\hfill\openbox

\medskip

In the next example we give a semigroup $S$ in which the identity relation $\iota _S$ is a left equalizer simple congruence on $S$ such that the
right colon of ${\mathbb F}[\iota _S]$ with respect to ${\mathbb F}[S]$ equals the augmentation ideal of ${\mathbb F}[S]$.

\begin{example}\label{example3}\rm
Let $S=\{a, b, c, d\}$ be a semigroup defined by Table~\ref{Fig2}.

\medskip

\begin{table}[htbp]
\begin{center}
\begin{tabular}{l|l l l l }
 &$a$&$b$&$c$&$d$\\ \hline
$a$&$a$&$a$&$a$&$a$\\
$b$&$b$&$b$&$b$&$b$\\
$c$&$a$&$a$&$a$&$a$\\
$d$&$b$&$b$&$b$&$b$\\
\end{tabular}
\caption{}\label{Fig2}
\end{center}
\end{table}

It is easy to see that $S$ can be obtained by applying Construction 1 in that case when $T=\{ e\}$ is a one-element semigroup, $S_e=\{ a, b, c, d\}$
and $f_{(e, e)}$ defined by
\[(a)f_{(e, e)}=a,\quad (b)f_{(e, e)}=b,\quad (c)f_{(e, e)}=a,\quad (d)f_{(e, e)}=b.\] It is a matter of checking to see that
\[f_{(e, e)}\circ f_{(e, e)}=f_{(e, e)}\] and so the conditions of Construction 1 are satisfied. The semigroup $S$ is left equalizer simple and so the
identity relation $\iota _S$ of $S$ is a left equalizer simple congruence on $S$.
For an element $\alpha a+\beta b+\gamma c+\delta d\in {\mathbb F}[S]$,
\[\alpha a+\beta b+\gamma c+\delta d\in {\mathbb F}[\iota ^{(1)}_S]=({\mathbb F}[\iota _S]:_r{\mathbb F}[S])\]
if and only if \[(\alpha +\beta +\gamma +\delta)=0.\] Thus
the right colon $({\mathbb F}[\iota _S]:_r{\mathbb F}[S])$ equals the augmentation ideal ${\mathbb F}[\omega _S]$.
It is easy to see that
$S$ is an ideal extension of the left zero semigroup $L=\{a, b\}$ by the null semigroup $S/L$.
\end{example}

\medskip

{\bf Acknowledgement}: I would like to thank the referee for the valuable remarks on the original version of the paper.

\medskip




\begin{thebibliography}{1}

\bibitem{Clifford1:sg-1} A.H. Clifford and G.B. Preston, {\it The Algebraic Theory of Semigroups I}, American Mathematical Society, Providence, R. I., 1961
\bibitem{Bogdanovic:sg-2} S. Bogdanovi\v c, {\it Semigroups with a system of subsemigroups}, University of Novi Sad, 1985
\bibitem{Cr1:sg-3} J. L. Chrislock, {\it Semigroups Whose Regular Representation is a Group}, Proc. Japan Acad., 40(1964), no. 10, 799-800
\bibitem{Cr2:sg-4} J. L. Chrislock, {\it Semigroups Whose Regular Representation is a Right Group}, The American Mathematical Monthly, 74(1967), No. 9, 1097-1100
\bibitem{Nagybook:sg-5} A. Nagy, {\it Special Classes of Semigroups}, Kluwer Academic Publishers, Dordrecht, Boston, London, 2001
\bibitem{JespersOkn:sg-12}, E. Jespers and J. Okni\'nski, {\it Noetherian Semigroup Algebras}, Springer-Verlag, Dordrecht, 2007
\bibitem{Nagyleft:sg-6} A. Nagy, {\it Left reductive congruences on semigroups}, Semigroup Forum, 87(2013), 129-148
\bibitem{Nagy0:sg-7} {\it On faithful representations of finite semigroups ${S}$ of degree $|{S}|$ over the fields}, International Journal of Algebra, 7(2013), 115-129
\bibitem{NagyRonyai:sg-8} A. Nagy and L. R\'onyai, {\it Finite Semigroups whose Semigroup Algebra over a Field Has a Trivial Right Annihilator},
International Journal of Contemporary Mathematical Science, 9(2014), 25-36
\bibitem{Okninski:sg-9} J. Okni\'nski, {\it Semigroup Algebras}, Marcel Dekker, Inc., New York, 1991
\bibitem{Y:sg-10} M. Yamada, {\it A not on middle unitary semigroups}, K\"odai Math. Sem. Rep., 7(1955), 49-52
\bibitem{Howie:sg-11} J.M. Howie, {\it An Introduction to Semigroup Theory}, Academic Press, London, 1976
\end{thebibliography}
\end{document}